\documentclass{fundam}

\usepackage[latin1]{inputenc}
\usepackage{amssymb}

\newtheorem{The}{Theorem}[section]
\newtheorem{Pro}[The]{Proposition}
\newtheorem{Deff}[The]{Definition}
\newtheorem{Lem}[The]{Lemma}
\newtheorem{Rem}[The]{Remark}
\newtheorem{Exa}[The]{Example}
\newtheorem{Cor}[The]{Corollary}

\newcommand{\fa}{\forall}

\newcommand{\Si}{\Sigma}
\newcommand{\Sis}{\Sigma^\star}
\newcommand{\Sio}{\Sigma^\omega}
\newcommand{\ra}{\rightarrow}
\newcommand{\hs}{\hspace{12mm}

\noi}
\newcommand{\lra}{\leftrightarrow}

\newcommand{\ite}{\item}

\newcommand{\ol}{ $\omega$-language}

\newcommand{\om}{\omega}
\newcommand{\nl}{\newline}
\newcommand{\noi}{\noindent}

\newcommand{\proo}{\noi {\bf Proof.} }
\newcommand {\ep}{\hfill $\square$}

\begin{document}

\setcounter{page}{1}
\issue{(2009)}

\title{{\bf  On  Recognizable Tree Languages  \\   Beyond the Borel Hierarchy }}

\author{Olivier Finkel \\ {\it Equipe de Logique Math\'ematique}  
  \\  CNRS et Universit\'e Paris 7,  France.\\ finkel@logique.jussieu.fr
\and 
  Pierre Simonnet \\ {\it   Systèmes physiques pour l'environnement }  \\
 Facult\'e des Sciences, Universit\'e de Corse\\
 {\it Quartier Grossetti BP52 20250, Corte, France } 
\\ simonnet@univ-corse.fr}

\maketitle

\runninghead{Olivier Finkel, Pierre Simonnet }{On  Recognizable Tree Languages   Beyond the Borel Hierarchy}

\begin{abstract}
\noi We investigate the topological complexity of  non Borel recognizable tree languages with regard to the difference hierarchy of analytic sets. 
We show that, for each integer $n \geq 1$, there is a $D_{\om^n}({\bf \Si}^1_1)$-complete  tree  language $\mathcal{L}_n$  
accepted by a (non deterministic)  Muller tree automaton. 
On the other hand, we prove that a tree language accepted by an unambiguous B\"uchi tree automaton must be Borel. 
Then we consider the game tree languages $W_{(\iota, \kappa)}$, for Mostowski-Rabin indices $(\iota, \kappa)$. We prove that 
the $D_{\om^n}({\bf \Si}^1_1)$-complete tree languages $\mathcal{L}_n$ are Wadge reducible to the game tree language
$W_{(\iota, \kappa)}$ for $\kappa-\iota \geq 2$. 
In particular these languages  $W_{(\iota, \kappa)}$ are not in any class $D_{\alpha}({\bf \Si}^1_1)$ 
for $\alpha < \om^\om$. 
\end{abstract}


\section{Introduction}

\noi   A way to study the complexity of languages of infinite words or infinite trees accepted by various kinds of automata is to study their topological 
complexity, and firstly to locate them with regard to the Borel and the projective hierarchies. 
 It is well known that every $\om$-language accepted by a deterministic B\"uchi automaton is a ${\bf \Pi}^0_2$-set. This implies  that 
  any $\om$-language accepted by a deterministic   Muller automaton is a 
boolean combination of ${\bf \Pi}^0_2$-sets hence a
${\bf \Delta}^0_3$-set.      \cite{Thomas90,Staiger97,PerrinPin}. But then it follows from Mc Naughton's Theorem, that 
all regular $\om$-languages, which are accepted by deterministic Muller automata, are also 
${\bf \Delta}^0_3$-sets. 
The Borel hierarchy of regular $\om$-languages is then determined. Moreover  Wagner  determined a much more refined hierarchy on regular 
$\om$-languages, which is in fact the 
trace of the Wadge hierarchy on  regular $\om$-languages, now called the Wagner hierarchy. 

\hs On the other hand,  many questions remain open about the topological 
complexity of regular languages of infinite trees. We know that they can be much more complex than regular sets of infinite words. 
Skurczynski proved that for every integer $n\geq 1$, there are some   ${\bf \Pi}^0_n$-complete and some       ${\bf \Si}^0_n$-complete  
regular tree languages,            \cite{Skurczynski93}.  Notice that it is an open question to know whether there exist some regular sets of trees 
which are Borel sets of infinite rank.   But there exist some regular sets of trees which are not Borel. Niwinski showed that there are some 
${\bf \Si}^1_1$-complete regular sets of trees accepted by B\"uchi tree automata, and some 
${\bf \Pi}^1_1$-complete regular sets of trees accepted by deterministic Muller  tree automata, 
 \cite{Niwinski85}. Every set of trees accepted by a  B\"uchi tree
automaton is a ${\bf \Si}^1_1$-set and every  set of trees accepted by a deterministic Muller  tree automaton is a 
${\bf \Pi}^1_1$-set. Niwinski and Walukiewicz proved that a tree language which is accepted by a deterministic Muller  tree automaton 
is either in the class ${\bf \Pi}^0_3$ or ${\bf \Pi}^1_1$-complete,    \cite{NiwiskiWalukiewicz03}.  More recent  results of Duparc and Murlak, 
 on the Wadge hierarchy of recognizable  tree languages, may be found in \cite{Murlak-LMCS, ADMN}.

\hs      It follows from the definition of acceptance by  non deterministic Muller or Rabin automata and from Rabin's complementation Theorem that 
every regular set of trees is a ${\bf \Delta}^1_2$-set, see \cite{Rabin69,PerrinPin,Thomas90,LescowThomas}. 
But there are only few known results on the complexity of non Borel regular tree languages. 
The second author gave examples of   $D_{\om^n}({\bf \Si}^1_1)$-complete regular  tree  languages in \cite{Simonnet92}. 
Arnold and Niwinski showed in \cite{ArnoldNiwinski08} that the game tree languages  $W_{(\iota, \kappa)}$ form a infinite hierarchy of non Borel 
regular sets of trees with regard to the Wadge 
reducibility. 

\hs In this paper, we investigate the topological complexity of  non Borel recognizable tree languages with regard to the difference hierarchy of analytic sets. 
We show that, for each integer $n \geq 1$, there is a $D_{\om^n}({\bf \Si}^1_1)$-complete  tree  language $\mathcal{L}_n$  
accepted by a (non deterministic)  Muller tree automaton. 
On the other hand, we prove that non Borel recognizable tree languages accepted by B\"uchi tree automata have the maximum degree of ambiguity. 
In particular, a tree language recognized by an unambiguous B\"uchi tree automaton must be Borel. 
Then we consider the game tree languages $W_{(\iota, \kappa)}$, for Mostowski-Rabin indices $(\iota, \kappa)$. We prove that 
the $D_{\om^n}({\bf \Si}^1_1)$-complete tree languages $\mathcal{L}_n$ are Wadge reducible to the game tree language
$W_{(\iota, \kappa)}$ for $\kappa-\iota \geq 2$. 
In particular,  these languages  $W_{(\iota, \kappa)}$ are not in any class $D_{\alpha}({\bf \Si}^1_1)$ 
for $\alpha < \om^\om$. 

\hs The paper is organized as follows. In Section 2 we recall the notions of B\"uchi or Muller tree automata and of regular tree languages. The notions of topology,
including the definition of the difference hierarchy of analytic sets, are recalled in Section 3.  We show in Section  4 that there are 
 $D_{\om^n}({\bf \Si}^1_1)$-complete  tree  languages $\mathcal{L}_n$  
accepted by   Muller tree automata. We consider the complexity of game tree languages in Section 5.

\section{Recognizable tree languages}

We  recall now  usual notations of formal language theory. 
\nl  When $\Si$ is a finite alphabet, a {\it non-empty finite word} over $\Si$ is any 
sequence $x=a_1\cdots a_k$, where $a_i\in\Sigma$ 
for $i=1,\ldots ,k$ , and  $k$ is an integer $\geq 1$. The {\it length}
 of $x$ is $k$, denoted by $|x|$.
 The {\it empty word} has no letter and is denoted by $\lambda$; its length is $0$. 
 $\Sis$  is the {\it set of finite words} (including the empty word) over $\Sigma$. A {\it   finitary language}   
over the alphabet $\Si$   is a subset of 
 $\Sis$. 

 \hs  The {\it first infinite ordinal} is $\om$.
 An $\om$-{\it word} over $\Si$ is an $\om$ -sequence $a_1 \cdots a_n \cdots$, where for all 
integers $ i\geq 1$, ~
$a_i \in\Sigma$.  When $\sigma$ is an $\om$-word over $\Si$, we write
 $\sigma =\sigma(1)\sigma(2)\cdots \sigma(n) \cdots $,  where for all $i$,~ $\sigma(i)\in \Si$,
and $\sigma[n]=\sigma(1)\sigma(2)\cdots \sigma(n)$  for all $n\geq 1$ and $\sigma[0]=\lambda$.
\nl 
 The usual concatenation product of two finite words $u$ and $v$ is 
denoted $u\cdot v$ (and sometimes just $uv$). This product is extended to the product of a 
finite word $u$ and an $\om$-word $v$: the infinite word $u\cdot v$ is then the $\om$-word such that:
\nl $(u\cdot v)(k)=u(k)$  if $k\leq |u|$ , and 
 $(u\cdot v)(k)=v(k-|u|)$  if $k>|u|$.
\nl   The {\it prefix relation} is denoted $\sqsubseteq$: a finite word $u$ is a {\it prefix} 
of a finite word $v$ (respectively,  an infinite word $v$), denoted $u\sqsubseteq v$,  
 if and only if there exists a finite word $w$ 
(respectively,  an infinite word $w$), such that $v=u\cdot w$.
\nl   
 The {\it set of } $\om$-{\it words} over  the alphabet $\Si$ is denoted by $\Si^\om$.
An  $\om$-{\it language} over an alphabet $\Sigma$ is a subset of  $\Si^\om$.

\hs We introduce now  languages of infinite binary trees whose nodes
are labelled in a finite alphabet $\Si$.
\nl A node of an infinite binary tree is represented by a finite  word over 
the alphabet $\{l, r\}$ where $r$ means ``right" and $l$ means ``left". Then an 
infinite binary tree whose nodes are labelled  in $\Si$ is identified with a function
$t: \{l, r\}^\star \ra \Si$. The set of  infinite binary trees labelled in $\Si$ will be 
denoted $T_\Si^\om$.

\hs Let $t$ be a tree. A branch $B$ of $t$ is a subset of the set of nodes of $t$ which 
is linearly ordered by the tree partial order $\sqsubseteq$ and which 
is closed under prefix relation, 
i.e. if  $x$ and $y$ are nodes of $t$ such that $y\in B$ and $x \sqsubseteq y$ then $x\in B$.
\nl A branch $B$ of a tree is said to be maximal iff there is not any other branch of $t$ 
which strictly contains $B$.

\hs     Let $t$ be an infinite binary tree in $T_\Si^\om$. If $B$ is a maximal branch of $t$,
then this branch is infinite. Let $(u_i)_{i\geq 0}$ be the enumeration of the nodes in $B$
which is strictly increasing for the prefix order. 
\nl  The infinite sequence of labels of the nodes of  such a maximal 
branch $B$, i.e. $t(u_0)t(u_1) \cdots t(u_n) \cdots $  is called a path. It is an $\om$-word 
over the alphabet $\Si$.

\hs     Let then $L\subseteq \Si^\om$ be an \ol~ over $\Si$. Then  we denote $\exists \mathrm{Path}(L)$  the set of 
infinite trees $t$ in $T_\Si^\om$ such that $t$ has (at least) one path in $L$.

\hs We are now going to define tree automata and recognizable tree languages. 

\begin{Deff} A (nondeterministic topdown) tree automaton  is a quadruple $\mathcal{A}=(K,\Si,\Delta, q_0)$, where $K$ 
is a finite set of states, $\Sigma$ is a finite input alphabet, $q_0 \in K$ is the initial state
and $\Delta \subseteq  K \times   \Si   \times  K \times   K$ is the transition relation. The tree automaton  $\mathcal{A}$ is said to be 
deterministic if the relation  $\Delta$ is a functional one, i.e.  if  for each $(q, a) \in K \times   \Si $ there is at most one pair of states $(q', q'')$ 
such that  $(q, a, q', q'') \in \Delta$. 
\nl A run of the tree automaton  $\mathcal{A}$ on an infinite binary tree $t\in T_\Si^\om$ is a infinite binary tree $\rho \in T_K^\om$ such that: 
\nl (a)  $\rho (\lambda)=q_0$  and  ~~(b) for each $u \in \{l, r\}^\star$,  $(\rho(u), t(u), \rho(u.l), \rho(u.r))  \in \Delta$. 
\end{Deff}

\begin{Deff}
A B\"uchi (nondeterministic topdown) tree automaton  is a  5-tuple $\mathcal{A}=(K,\Si,\Delta, q_0, F)$, where $(K,\Si,\Delta, q_0)$ is a 
tree automaton and $F \subseteq K$ is the set of accepting states. 
\nl A run $\rho$ of the  B\"uchi tree automaton $\mathcal{A}$ on an infinite binary tree $t\in T_\Si^\om$ is said to be accepting if  
for each path of $\rho$ there is some accepting state appearing infinitely  often on this path. 
\nl The tree language $L(\mathcal{A})$ accepted by the B\"uchi tree automaton $\mathcal{A}$ is the set of infinite binary trees $t\in T_\Si^\om$ 
such that there is (at least) one accepting run of $\mathcal{A}$ on $t$. 
\end{Deff}

\begin{Deff}
A Muller  (nondeterministic topdown) tree automaton  is a  5-tuple $\mathcal{A}=(K,\Si,\Delta, q_0, \mathcal{F})$, where $(K,\Si,\Delta, q_0)$ is a 
tree automaton and $\mathcal{F} \subseteq 2^K$ is the collection of  designated state sets. 
\nl A run $\rho$ of the  Muller  tree automaton $\mathcal{A}$ on an infinite binary tree $t\in T_\Si^\om$ is said to be accepting if  
for each path $p$ of $\rho$, the set of   states appearing infinitely  often on this path is in $\mathcal{F}$. 
\nl The tree language $L(\mathcal{A})$ accepted by the  Muller tree automaton $\mathcal{A}$ is the set of infinite binary trees $t\in T_\Si^\om$ 
such that there is (at least) one accepting run of $\mathcal{A}$ on $t$. 
\nl The class {\bf $REG$} of regular, or recognizable, tree languages is the class of  tree  languages accepted by some Muller automaton. 
\end{Deff}

\begin{Rem}
Each tree language accepted by some (deterministic)  B\"uchi automaton is also accepted by some (deterministic) 
Muller automaton. A tree language is accepted by a Muller tree  automaton iff it is accepted by some Rabin tree automaton. We refer for instance to 
\cite{Thomas90,PerrinPin} for the definition of Rabin tree  automaton. 
\end{Rem}

\begin{Exa}
Let  $L\subseteq \Si^\om$ be a regular $\om$-language 
(see \cite{PerrinPin} about regular $\om$-languages which are the $\om$-languages accepted by  B\"uchi or Muller automata). 
Then  the set  $\exists \mathrm{Path}(L) \subseteq  T_\Si^\om$ is accepted by a 
B\"uchi tree automaton, hence also by a Muller  tree automaton. 
\nl The set of infinite binary trees $t\in T_\Si^\om$ having all their paths in $L$, denoted $\fa  \mathrm{Path}(L)$, 
  is accepted by a deterministic Muller  tree automaton. It is in fact the 
complement of the set $\exists \mathrm{Path}(\Sio-L) $. 
\end{Exa}

\section{Topology}

\hs We assume the reader to be familiar with basic notions of topology which
may be found in \cite{Moschovakis80,LescowThomas,Kechris94,Staiger97,PerrinPin}.
There is a natural metric on the set $\Sio$ of  infinite words 
over a finite alphabet 
$\Si$ containing at least two letters which is called the {\it prefix metric} and defined as follows. For $u, v \in \Sio$ and 
$u\neq v$ let $\delta(u, v)=2^{-l_{\mathrm{pref}(u,v)}}$ where $l_{\mathrm{pref}(u,v)}$ 
 is the first integer $n$
such that the $(n+1)^{st}$ letter of $u$ is different from the $(n+1)^{st}$ letter of $v$. 
This metric induces on $\Sio$ the usual  Cantor topology for which {\it open subsets} of 
$\Sio$ are in the form $W\cdot \Si^\om$, where $W\subseteq \Sis$.
A set $L\subseteq \Si^\om$ is a {\it closed set} iff its complement $\Si^\om - L$ 
is an open set.

 \hs     There is also a natural topology on the set $T_\Si^\om$ \cite{Moschovakis80,LescowThomas,Kechris94}. 
It is defined 
by the following distance. Let $t$ and $s$ be two distinct infinite trees in $T_\Si^\om$. 
Then the distance between $t$ and $s$ is $\frac{1}{2^n}$ where $n$ is the smallest integer 
such that $t(x)\neq s(x)$ for some word $x\in \{l, r\}^\star$ of length $n$.
\nl The open sets are then in the form $T_0\cdot T_\Si^\om$ where $T_0$ is a set of finite labelled
trees. $T_0\cdot T_\Si^\om$ is the set of infinite binary trees 
which extend some finite labelled binary tree $t_0\in T_0$, $t_0$ is here a sort of prefix, 
an ``initial subtree"
of a tree in $t_0\cdot T_\Si^\om$.
\nl It is well known that the set $T_\Si^\om$, equipped with this topology, is homeomorphic to the Cantor set $2^\om$, hence  also to the topological spaces 
$\Sio$, where $\Si$ is an alphabet having at least two letters. 

\hs We  now define the {\it Borel Hierarchy} of subsets of $\Si^\om$. It is defined similarly on the space $T_\Si^\om$.

\begin{Deff}
For a non-null countable ordinal $\alpha$, the classes ${\bf \Si}^0_\alpha$
 and ${\bf \Pi}^0_\alpha$ of the Borel Hierarchy on the topological space $\Si^\om$ 
are defined as follows:
\nl ${\bf \Si}^0_1$ is the class of open subsets of $\Si^\om$, 
 ${\bf \Pi}^0_1$ is the class of closed subsets of $\Si^\om$, 
\nl and for any countable ordinal $\alpha \geq 2$: 
\nl ${\bf \Si}^0_\alpha$ is the class of countable unions of subsets of $\Si^\om$ in 
$\bigcup_{\gamma <\alpha}{\bf \Pi}^0_\gamma$.
 \nl ${\bf \Pi}^0_\alpha$ is the class of countable intersections of subsets of $\Si^\om$ in 
$\bigcup_{\gamma <\alpha}{\bf \Si}^0_\gamma$.
\end{Deff}

\noi For 
a countable ordinal $\alpha$,  a subset of $\Si^\om$ is a Borel set of {\it rank} $\alpha$ iff 
it is in ${\bf \Si}^0_{\alpha}\cup {\bf \Pi}^0_{\alpha}$ but not in 
$\bigcup_{\gamma <\alpha}({\bf \Si}^0_\gamma \cup {\bf \Pi}^0_\gamma)$.

\hs    
There exists another hierarchy beyond the Borel hierarchy, which is called the 
projective hierarchy. The classes ${\bf \Si}^1_n$ and ${\bf \Pi}^1_n$, for integers $n\geq 1$, 
of the projective hierarchy are  obtained from  the Borel hierarchy by 
successive applications of operations of projection and complementation.
The first level of the projective hierarchy is formed by the class ${\bf \Si}^1_1$  of  {\it analytic sets} and the 
class ${\bf \Pi}^1_1$ of {\it co-analytic sets} which are complements of 
analytic sets.  
In particular, 
the class of Borel subsets of $\Si^\om$ is strictly included in 
the class  ${\bf \Si}^1_1$ of {\it analytic sets} which are 
obtained by projection of Borel sets. 

\begin{Deff} 
A subset $A$ of  $\Si^\om$ is in the class ${\bf \Si}^1_1$ of {\bf analytic} sets
iff there exists another finite set $Y$ and a Borel subset $B$  of  $(\Si \times Y)^\om$ 
such that $ x \in A \lra \exists y \in Y^\om $ such that $(x, y) \in B$, 
where $(x, y)$ is the infinite word over the alphabet $\Si \times Y$ such that
$(x, y)(i)=(x(i),y(i))$ for each  integer $i\geq 1$.
\end{Deff}

\begin{Rem}
In the above definition we could take $B$ in the class ${\bf \Pi}^0_2$. Moreover 
analytic subsets of $\Si^\om$ are the projections of  ${\bf \Pi}^0_1$-subsets of 
$\Si^\om \times \om^\om$, where $ \om^\om$ is the Baire space, \cite{Moschovakis80}. 
\end{Rem}

\noi  We now define the notion of Wadge reducibility  via the  reduction by continuous functions. 
Let $X$, $Y$ be two finite alphabets. 
For $L\subseteq X^\om$ and $L'\subseteq Y^\om$, $L$ is said to be Wadge reducible to $L'$, denoted by $L\leq _W L'$, 
 iff there exists a continuous function $f: X^\om \ra Y^\om$, such that
$L=f^{-1}(L')$.

\hs We now define completeness with regard to reduction by continuous functions. 
For a countable ordinal  $\alpha\geq 1$, and an integer $n\geq 1$,  a set $F\subseteq \Si^\om$ is said to be 
a ${\bf \Si}^0_\alpha$  
(respectively,  ${\bf \Pi}^0_\alpha$, ${\bf \Si}^1_n$, ${\bf \Pi}^1_n$)-{\it complete set} 
iff for any set $E\subseteq Y^\om$  (with $Y$ a finite alphabet): 
 $E\in {\bf \Si}^0_\alpha$ (respectively,  $E\in {\bf \Pi}^0_\alpha$,  $E\in {\bf \Si}^1_n$, $E\in {\bf \Pi}^1_n$) 
iff $E \leq_W F$. 
 ${\bf \Si}^0_n$
 (respectively ${\bf \Pi}^0_n$)-complete sets, with $n$ an integer $\geq 1$, 
 are thoroughly characterized in \cite{Staiger86a}.  

  \hs   The Borel hierarchy and the projective hierarchy on $T_\Si^\om$ are defined from open 
sets in the same manner as in the case of the topological space $\Si^\om$. 

\hs The $\om$-language  $\mathcal{R}=(0^\star\cdot 1)^\om$  
is a well known example of 
${\bf \Pi}^0_2 $-complete subset of $\{0, 1\}^\om$. It is the set of 
$\om$-words over $\{0, 1\}$ having infinitely many occurrences of the letter $1$. 
Its  complement 
$\{0, 1\}^\om - (0^\star\cdot 1)^\om$ is a 
${\bf \Si}^0_2 $-complete subset of $\{0, 1\}^\om$.
\nl The set of infinite trees in $T_\Si^\om$, where $\Si=\{0, 1\}$, having at least one  path in the $\om$-language $\mathcal{R}=(0^\star\cdot 1)^\om$ is 
${\bf \Si}^1_1$-complete. Its complement is the set of trees in $T_\Si^\om$ having all their paths in $\{0, 1\}^\om - (0^\star\cdot 1)^\om$; it is 
${\bf \Pi}^1_1$-complete. 

\hs We now recall the notion of  difference hierarchy of analytic sets. 
 Let $\eta\! <\!\omega_1$ (where $\omega_1$ is the first uncountable ordinal) 
be an ordinal and $(A_\theta )_{\theta <\eta}$ be an increasing sequence of subsets of some space $X$,
 then the  set 
$D_\eta [(A_\theta )_{\theta <\eta}]$ is the set of elements $x \in X$ such that $x\!\in\! A_\theta\!\setminus\!\bigcup_{\theta'<\theta}\ A_{\theta'}$ for some 
$\theta\! <\!\eta$ whose  parity  is opposite to that of $\eta$. (Recall that a countable ordinal $\gamma$ is said to be even iff it can be written in the form 
$\gamma=\alpha + n$, where $\alpha$ is a limit ordinal and $n$ is an even non-negative  integer; otherwise the ordinal $\gamma$ is said to be odd; notice that all 
limit ordinals, like the ordinals $\om^n$, $n \geq 1$, or $\om^\om$, are even ordinals.)

\hs 
 We can now define the class of $\eta$-differences of analytic subsets of $X$, where $X=\Si^\om$ or $X=T_\Si^\om$. 

\hs $D_\eta ({\bf \Si}^1_1 )\! :=\!\{ D_\eta 
[(A_\theta )_{\theta <\eta}]\mid \mbox{ for each ordinal } \theta <\eta  ~~ A_\theta  \mbox{ is a } {\bf \Si}^1_1\mbox{-set }\}$

\hs It is well known that the hierarchy of differences of  analytic sets is strict, i.e. that for all countable ordinals $\alpha < \beta < \omega_1$, it holds that 
$D_\alpha ({\bf \Si}^1_1)  \subset   D_\beta  ({\bf \Si}^1_1 )$. This is considered as a folklore result of descriptive set theory which follows from the 
existence of universal sets for each class $D_\alpha ({\bf \Si}^1_1)$. Indeed we know first that the class ${\bf \Si}^1_1$ of analytic sets admits a universal 
set, see \cite[page 205]{Kechris94}or \cite[page 43]{Moschovakis80}. 
Then, using classical methods of descriptive set theory, one can show that, for each countable ordinal $\alpha$, 
the class $D_\alpha ({\bf \Si}^1_1)$ admits also  a universal set,  see \cite[page 443]{Kanamori}. 
This implies, as in the case of the Borel hierarchy in \cite[page 168]{Kechris94}, 
that the difference hierarchy of  analytic sets is strict. As a universal set for the class $D_\alpha ({\bf \Si}^1_1)$ is also a 
$D_\alpha ({\bf \Si}^1_1)$-complete set  for reduction 
by continuous functions, this implies also that there exists a $D_\alpha ({\bf \Si}^1_1)$-complete set. 
\nl Notice that in the sequel we shall  only  consider the classes $D_\alpha ({\bf \Si}^1_1)$, for  ordinals $\alpha < \om^\om$, and that we shall 
reprove that there exists some $D_\alpha ({\bf \Si}^1_1)$-complete subsets of $T_\Si^\om$, giving examples which are regular sets of trees. 

\hs Another folklore result of descriptive set theory is that  the union  $\bigcup_{\alpha <\omega_1} D_\alpha ({\bf \Si}^1_1)$ 
represents only a small part of the class ${\bf \Delta}_2^1$. It is quoted for instance in \cite{Steel82} or \cite[page 443]{Kanamori}.  (It is noticed 
in \cite{Steel82}  that the union $\bigcup_{\alpha <\omega_1} D_\alpha ({\bf \Si}^1_1)$  is strictly included in the class 
$\mathcal{A}({\bf \Pi}^1_1)$ which is the closure of the class ${\bf \Pi}^1_1$ under Souslin's operation. The class $\mathcal{A}({\bf \Pi}^1_1)$ is 
 included in the class ${\bf \Delta}_2^1$ by \cite[2.B.5  page 75]{Moschovakis80}).
Notice however that this result is not necessary 
in the sequel. 

\section{$D_{\om^n}({\bf \Si}^1_1)$-complete recognizable languages}

\noi    It follows from the definition of the B\"uchi acceptance condition for infinite trees that each tree language recognized by a (non deterministic) 
B\"uchi tree automaton is an analytic set. 
\nl Niwinski showed that some B\"uchi     recognized   tree languages  are actually  ${\bf \Si}^1_1$-complete sets.  
An example is any tree language $T \subseteq T_\Si^\om$  in the form $\exists \mathrm{Path}(L)$, 
where $L\subseteq \Si^\om$ is a regular $\om$-language 
 which is a 
${\bf \Pi}^0_2$-complete subset of $\Si^\om$. 
In particular,   the tree language $\mathcal{L}=\exists \mathrm{Path}( \mathcal{R} )$, where 
$\mathcal{R}=(0^\star\cdot 1)^\om$, is ${\bf \Si}^1_1$-complete 
hence non Borel      \cite{Niwinski85, PerrinPin,Simonnet92}.
\nl Notice that its complement   $\mathcal{L}^-=\fa \mathrm{Path}( \{0, 1\}^\om - (0^\star\cdot 1)^\om )$ is a ${\bf \Pi}^1_1$-complete set. 
It cannot be accepted by any 
B\"uchi tree automaton because it is not a ${\bf \Si}^1_1$-set. On the other hand, it can be easily seen that it is accepted by a deterministic 
Muller tree automaton. 

\hs The tree languages $\mathcal{L}$ and  $\mathcal{L}^-$ 
have been used by the second author in \cite{Simonnet92} to give examples of  
$D_{\om^n}({\bf \Si}^1_1)$-complete recognizable tree languages, for integers $n\geq 1$.  
We now give first the construction of a $D_{\om}({\bf \Si}^1_1)$-complete set. 

\hs For a tree $t \in  T_\Si^\om$ and $u\in \{l, r\}^\star$, we shall denote  $t_u : \{l, r\}^\star \ra \Si$ the subtree defined by 
$t_u(v)=t(u\cdot v)$ for all $v\in \{l, r\}^\star$. It is in fact the subtree of $t$ which is rooted in $u$. 

\hs Now we can define  a      $D_{\om}({\bf \Si}^1_1)$-complete    tree language  $\mathcal{L}_1$. 

\hs  $\mathcal{L}_1=\{ t \in T_{\{0, 1\}}^\om  \mid  \exists n\geq 0 ~~t_{l^n\cdot r} \in \mathcal{L} \mbox{ and  min} 
\{ n\geq 0 \mid t_{l^n\cdot r} \in \mathcal{L} \} 
\mbox{ is odd }\}$. 

\begin{Pro}\label{L1}
The tree language $\mathcal{L}_1$ is $D_{\om}({\bf \Si}^1_1)$-complete. 
\end{Pro}

\proo  We first show that the language $\mathcal{L}_1$  is in the class $D_{\om}({\bf \Si}^1_1)$. 
\nl Consider firstly, for some integer $k\geq0$,  the set $T_k=\{ t \in T_{\{0, 1\}}^\om  \mid   ~t_{l^{k}\cdot r} \in \mathcal{L} \}$.  
It is clear that this set is in the class 
${\bf \Si}^1_1$ because the function $F_k : T_{\{0, 1\}}^\om \ra T_{\{0, 1\}}^\om$ defined by $F_k(t)=t_{l^{k}\cdot r}$ is continuous 
and $T_k=F_k^{-1}(\mathcal{L})$ and the class ${\bf \Si}^1_1$ is closed under inverses of continuous functions. 
\nl Let now $H_n=\{ t \in T_{\{0, 1\}}^\om  \mid  \exists k \leq n ~t_{l^{k}\cdot r} \in \mathcal{L} \}$. This set is also in the class 
${\bf \Si}^1_1$ because the class ${\bf \Si}^1_1$  is closed under finite (and even countable) union and $H_n=\bigcup_{k \leq n} T_k$. 
\nl The sets $H_n$ form an increasing sequence of ${\bf \Si}^1_1$-sets, and we can check that 
$$\mathcal{L}_1=D_\om [(H_n )_{n <\om}]$$

\noi We now prove that $\mathcal{L}_1$ is $D_{\om}({\bf \Si}^1_1)$-{\it complete}. 
\nl Let $L \subseteq \Sio$ be a $D_{\om}({\bf \Si}^1_1)$-subset of $\Sio$, where $\Si$ is an alphabet having at least two letters. Then 
there is an increasing sequence $(A_n)_{n\in\om}$ of  ${\bf \Si}^1_1$-subsets of $\Sio$ such that $L=D_\om [(A_n )_{n <\om}]$. 
On the other hand, we know that the tree language $\mathcal{L}$ is ${\bf \Si}^1_1$-complete. Thus for each integer $n \geq 0$ there exists a continuous 
function $f_n: \Sio  \ra T_{\{0, 1\}}^\om$ such that $A_n = f_n^{-1}(\mathcal{L})$. 
\nl We now define a function $F:  \Sio  \ra T_{\{0, 1\}}^\om$ by : for all $x \in \Sio$,  for all integers $k\geq 0$, $F(x)(l^k)=0$ 
and $F(x)_{l^k\cdot r}=f_k(x)$. 
It is clear that the function $F$ is continuous because each function $f_k$ is continuous. 
\nl We can now check that for every $x \in \Sio$,  $x$ is in the set $L=D_\om [(A_n )_{n <\om}]$ iff there is an odd integer $n$ such that 
$x\in A_n \ \setminus \ \bigcup_{k<n}\ A_{k}$ iff  there is an odd integer $n$ such that $f_n(x)\in \mathcal{L}$ and for all $k<n$ ~
$f_k(x) \in \mathcal{L}^-$. 
\nl This means that    $x \in L=D_\om [(A_n )_{n <\om}]$     iff $F(x) \in \mathcal{L}_1$. 
\nl Finally we have shown, using the reduction $F$, that $L=D_\om [(A_n )_{n <\om}] \leq_W   \mathcal{L}_1$ and so 
the tree language $ \mathcal{L}_1$ is $D_{\om}({\bf \Si}^1_1)$-complete. 
\ep

\hs We can now generalize this construction to obtain some $D_{\om^n}({\bf \Si}^1_1)$-complete tree languages, for every integer $n\geq 1$. 
\nl Recall first that an ordinal $\alpha$ is strictly smaller than the ordinal $\om^n$, where $n\geq 2$ is an integer, if and only if it admits a Cantor Normal Form 
$$\alpha= \om^{n-1}\cdot a_{n-1} + \om^{n-2}\cdot a_{n-2} + \ldots + \om\cdot a_1 + a_0$$
\noi where $a_{n-1}, a_{n-2}, \ldots , a_0,$ are non-negative integers. In that case we shall denote 
$\mathrm{Ord}(a_{n-1}, a_{n-2}, \ldots , a_0)=\om^{n-1}\cdot a_{n-1} + \om^{n-2}\cdot a_{n-2} + \ldots + \om\cdot a_1 + a_0$.  
\nl Recall also that if $\alpha= \mathrm{Ord}(a_{n-1}, a_{n-2}, \ldots , a_0)$ and $\beta=\mathrm{Ord}(b_{n-1}, b_{n-2}, \ldots , b_0)$, 
then  $\alpha < \beta$ if and only 
if there is an integer $k$ such that $0 \leq k \leq n-1$ and $a_j=b_j$ for $n-1\geq j >k$ and $a_k < b_k$. 

\hs  We now define the tree language $\mathcal{L}_n$, for $n\geq 2$, as the set of trees  $t \in T_{\{0, 1\}}^\om$ for which 
there exist some integers $a_{n-1}, a_{n-2}, \ldots , a_0 \geq 0$ such that: 

\begin{enumerate}
\ite  $t_{l^{a_{n-1}}\cdot r\cdot l^{a_{n-2}}\cdot r \cdots l^{a_0}\cdot r}$ is in  $\mathcal{L}$
and the parity of  $\mathrm{Ord}(a_{n-1}, a_{n-2}, \ldots , a_0)$ is odd, 
\ite If $\mathrm{Ord}(b_{n-1}, b_{n-2}, \ldots , b_0) < \mathrm{Ord}(a_{n-1}, a_{n-2}, \ldots , a_0)$ then the tree  
\nl  $t_{l^{b_{n-1}}\cdot r\cdot l^{b_{n-2}}.r \cdots l^{b_0}\cdot r}$ is not in  $\mathcal{L}$. 
\end{enumerate}

\begin{Pro}
For each integer $n\geq 2$, the tree language $\mathcal{L}_n$ is $D_{\om^n}({\bf \Si}^1_1)$-complete. 
\end{Pro}

\proo The proof is a simple generalization of the proof of Proposition \ref{L1}. Notice that  we have to use 
the closure of the class ${\bf \Si}^1_1$  under  countable (and not only under finite)  union. Details are here left to the reader. 
\ep 

\hs The tree languages $\mathcal{L}_n$  can not be accepted by any B\"uchi tree automaton
because  each tree language accepted by a (non deterministic) 
B\"uchi tree automaton is an analytic set and $D_{\om^n}({\bf \Si}^1_1)$-complete sets, for $n\geq 1$, are not in the class ${\bf \Si}^1_1$. 
We are going to see that the tree languages $\mathcal{L}_n$ 
are accepted by Muller tree automata. 

\hs We now recall the following result proved by Niwinski in \cite{Niwinski85}, see also for instance \cite{PerrinPin,Thomas90}. 

\begin{Lem}\label{det}
The language $\mathcal{L}^-=\fa \mathrm{Path}( \{0, 1\}^\om - (0^\star.1)^\om )$ is a ${\bf \Pi}^1_1$-complete set 
accepted by a deterministic Muller tree automaton. 
\end{Lem}

\noi On the other hand, the tree language $\mathcal{L}$ is a  ${\bf \Si}^1_1$-complete set. Thus it is not a ${\bf \Pi}^1_1$-set otherwise it would be 
in the class ${\bf \Delta}^1_1={\bf \Si}^1_1 \cap {\bf \Pi}^1_1$ which is the class of Borel sets by Suslin's Theorem. 
But every tree language which is recognizable by a {\it deterministic} Muller tree automaton is a ${\bf \Pi}^1_1$-set therefore 
the tree language $\mathcal{L}$ can not be accepted by any {\it deterministic} Muller tree automaton. 
 However we can now state the following result.

\begin{Lem}
The language $\mathcal{L}$ is a  ${\bf \Si}^1_1$-complete set accepted by a non deterministic B\"uchi tree automaton, 
hence also by a non deterministic Muller tree automaton. 
\end{Lem}

\proo We recall informally how we can define a non-deterministic B\"uchi tree automaton $\mathcal{A}$ accepting the language $\mathcal{L}$. 
When reading a tree $t\in \mathcal{L}$, the automaton $\mathcal{A}$, using the non determinism, guesses an infinite branch of the tree. Then the 
automaton checks, using the B\"uchi acceptance condition, that the sequence of labels of nodes on this branch forms an $\om$-word in $(0^\star.1)^\om$, 
i.e. contains an infinite number of letters $1$. 
\ep 

\begin{Lem}
For each integer $n \geq 1$, the language $\mathcal{L}_n$ is   
 accepted by a (non deterministic)  Muller tree automaton. 
\end{Lem}

\proo We first  construct a non deterministic  Muller tree automaton $\mathcal{A}_1$  accepting 
 the language $\mathcal{L}_1$. 

\hs Recall that, for each  tree $t\in \mathcal{L}_1$, there exists a least integer $n\geq 0$ such that $t_{l^n\cdot r} \in \mathcal{L}$. 
This (odd) integer is defined {\it in a unique way}. One can now construct, from   
Muller tree automata $\mathcal{A}^-$ and $\mathcal{A}^+$ accepting the tree languages $\mathcal{L}^-$ and $\mathcal{L}$, a  
 Muller tree automaton $\mathcal{A}_1$ accepting the tree language $\mathcal{L}_1$. 
Using the non-determinism, the automaton $\mathcal{A}_1$ will guess the  (odd) integer $n\geq 0$ and then, 
using the behaviour of $\mathcal{A}^-$ and $\mathcal{A}^+$,  it will check that $t_{l^n\cdot r} \in \mathcal{L}$ and that,  for every integer 
$k<n$, $t_{l^k\cdot r} \notin \mathcal{L}$. 

\hs We now give the exact construction of the  non deterministic  Muller tree automaton $\mathcal{A}_1$. 
\nl Let $\Si=\{0, 1\}$ and $\mathcal{A}^-=(K,\Si,\Delta, q_0, \mathcal{F})$ be a (deterministic) Muller tree automaton accepting 
the tree language $\mathcal{L}^-$. 
\nl And let 
$\mathcal{A}^+=(K',\Si,\Delta', q'_0, \mathcal{F}')$ be a (non deterministic) Muller tree automaton accepting 
the tree language $\mathcal{L}$. We assume that $K \cap K' =\emptyset$. 
\nl Then it is easy to see that the tree language $\mathcal{L}_1$ is accepted by the Muller tree automaton 
$\mathcal{A}_1=(K^1,\Si,\Delta^1,$ $ q^1_0, \mathcal{F}^1)$, where 
\nl $K^1=K \cup K'  \cup \{ q^1_0, q^1_1, q_f \}$, 
\nl $\Delta^1=\Delta \cup \Delta'  \cup \{ (q^1_0, a, q^1_1, q_0), (q^1_1, a, q_f, q'_0), (q_f, a, q_f, q_f), (q^1_1, a, q^1_0, q_0) \mid a\in \{0, 1\} \}$, 
\nl $\mathcal{F}^1=\mathcal{F}\cup \mathcal{F}' \cup \{q_f \}$.

\hs For every integer $n>1$, we can construct in a similar way a Muller tree automaton  $\mathcal{A}_n$ accepting the tree language 
$\mathcal{L}_n$. 

\hs Recall that for each  tree $t\in \mathcal{L}_n$ there exists a least ordinal $\alpha=\mathrm{Ord}(a_{n-1}, a_{n-2}, \ldots , a_0)<\om^n$ such that 
$t_{l^{a_{n-1}}\cdot r\cdot l^{a_{n-2}}\cdot r \cdots l^{a_0}\cdot r}$ is in  $\mathcal{L}$. This (odd) ordinal is defined {\it in a unique way}. 
\nl One can now construct, from the  
Muller tree automata $\mathcal{A}^-$ and $\mathcal{A}^+$ accepting the tree languages $\mathcal{L}^-$ and $\mathcal{L}$, a  
 Muller tree automaton $\mathcal{A}_n$ accepting the tree language $\mathcal{L}_n$. 
Using the non-determinism, the automaton $\mathcal{A}_n$ will guess the (odd) ordinal  $\alpha=\mathrm{Ord}(a_{n-1}, a_{n-2}, \ldots , a_0)<\om^n$ 
and then, using the behaviour of $\mathcal{A}^-$ and $\mathcal{A}^+$,  it will check that 
$t_{l^{a_{n-1}}\cdot r\cdot l^{a_{n-2}}\cdot r \cdots l^{a_0}\cdot r}$ is in  $\mathcal{L}$ 
and that for each ordinal $\beta=\mathrm{Ord}(b_{n-1}, b_{n-2}, \ldots , b_0) < \mathrm{Ord}(a_{n-1}, a_{n-2}, \ldots , a_0)$ the tree language 
$t_{l^{b_{n-1}}\cdot r\cdot l^{b_{n-2}}.r \cdots l^{b_0}\cdot r}$ is not in  $\mathcal{L}$. 

\ep 

\hs We can now summarize the above  results in the following theorem. 

\begin{The}\label{Ln}
For each integer $n \geq 1$, the language $\mathcal{L}_n$ is a    $D_{\om^n}({\bf \Si}^1_1)$-complete set
accepted by a (non deterministic) Muller tree automaton. 
\end{The}

\begin{Cor}
The class of tree languages recognized by  Muller tree automata is not included into the boolean closure of the class 
of tree languages recognized by B\"uchi tree automata. 
\end{Cor}

\proo We know that every tree language recognized by a B\"uchi tree automaton is a ${\bf \Si}^1_1$-set. But a tree language which is a boolean 
combination of ${\bf \Si}^1_1$-sets is in the class  $D_{\om}({\bf \Si}^1_1)$ which does not contain all 
tree languages recognized by  (non deterministic)  Muller tree automata. 
\ep 

\begin{Rem}{\rm 
We have given above  examples of    $D_{\om^n}({\bf \Si}^1_1)$-complete tree languages accepted by  Muller tree automata. 
In a similar way it is easy to   construct, for each ordinal $\alpha < \om^\om$,   a  $D_{\alpha}({\bf \Si}^1_1)$-complete tree language
accepted by a Muller tree automaton. 
Each ordinal  $\alpha < \om^\om$ may be written in the form $\alpha=\mathrm{Ord}(a_{n-1}, a_{n-2}, \ldots , a_0)<\om^n$ for some 
integer $n\geq 1$ and where $a_{n-1}, a_{n-2}, \ldots , a_0,$ are non-negative integers with $a_{n-1} \neq 0$. 
\nl The tree language $\mathcal{T}_\alpha$ is then the set of trees   $t \in T_{\{0, 1\}}^\om$ for which 
there exist some integers $b_{n-1}, b_{n-2}, \ldots , b_0 \geq 0$ such that: 

\begin{enumerate}
\ite $\mathrm{Ord}(b_{n-1}, b_{n-2}, \ldots , b_0) < \mathrm{Ord}(a_{n-1}, a_{n-2}, \ldots , a_0)$. 
\ite $t_{l^{b_{n-1}}\cdot r\cdot l^{b_{n-2}}.r \cdots l^{b_0}\cdot r}$ is in  $\mathcal{L}$
and the parity of  $\mathrm{Ord}(b_{n-1}, b_{n-2}, \ldots , b_0)$ is odd iff the parity of  
$\mathrm{Ord}(a_{n-1}, a_{n-2}, \ldots , a_0)$ is even. 
\ite If $\mathrm{Ord}(c_{n-1}, c_{n-2}, \ldots , c_0) < \mathrm{Ord}(b_{n-1}, b_{n-2}, \ldots , b_0)$ then the tree  
\nl  $t_{l^{c_{n-1}}\cdot r\cdot l^{c_{n-2}}.r \cdots l^{c_0}\cdot r}$ is not in  $\mathcal{L}$. 
\end{enumerate}

\noi The tree language $\mathcal{T}_\alpha$ is $D_{\alpha}({\bf \Si}^1_1)$-complete and  it is accepted by a (non deterministic)   
 Muller tree automaton. }

\end{Rem}

\noi  The above results show that the topological complexity of  tree languages recognized by  {\it non deterministic}
Muller tree automata is much greater than that of tree languages accepted by {\it deterministic} Muller tree automata. 

\hs Recall that a B\"uchi (respectively, Muller) tree automaton $\mathcal{A}$, reading trees labelled in the alphabet $\Si$, 
 is said to be unambiguous if and only if 
each tree $t \in T_\Si^\om$ admits at most one accepting run of $\mathcal{A}$. 

\hs  A natural question is whether the tree languages $\mathcal{L}_n$ could be accepted by unambiguous Muller  tree automata.  
A first step would be to prove that the 
tree language $\mathcal{L}$  is accepted by an unambiguous Muller  tree automaton.  
But this is not possible. We have learned    by 
personal communication from Damian Niwinski  that  the
language $\mathcal{L}$   is  inherently ambiguous, \cite{Niwinski09}. 

\hs We consider now  the notion of ambiguity for B\"uchi tree automata and we shall prove in particular that 
 a tree language accepted by an unambiguous B\"uchi tree automaton must be Borel. We shall indicate also why our methods do not 
work in the case of Muller automata. 

\hs 
We first  recall some notations and a  lemma proved in \cite{Fink-Sim}. 
\nl  For two finite alphabets  $\Si$ and $X$, 
if $B \subseteq \Sio  \times  X^\om $ and $\alpha 
\in \Sio$, we denote 
 $B_\alpha=\{\beta \in X^\om \mid (\alpha,\beta) \in B\}$ and   
$\mathrm{PROJ}_{\Sio}(B)=\{\alpha \in \Sio \mid B_\alpha \neq \emptyset \}$. 
\nl The cardinal of the continuum will be  denoted by  $2^{\aleph_0}$; it is also the cardinal 
of every set $\Sio$ or $T_\Si^\om$, where $\Si$ is an alphabet having at least two letters.

\begin{Lem}[\cite{Fink-Sim}]\label{lemma}
Let $\Si$ and $X$ be two finite alphabets having at least two letters and 
 $B$ be a Borel subset of 
$\Sio \times X^\om$ such that $\mathrm{PROJ}_{\Sio}(B)$ is not a Borel subset of $\Sio$.
Then there are $2^{\aleph_0}$ $\om$-words  $\alpha \in \Sio$ such that the section $B_\alpha$ 
has cardinality $2^{\aleph_0}$.
\end{Lem}

\proo Let $\Si$ and $X$ be two finite alphabets having at least two letters and 
 $B$ be a Borel subset of 
$\Sio \times X^\om$ such that $\mathrm{PROJ}_{\Sio}(B)$ is not  Borel.  

\hs In a first step we  prove that 
there are uncountably many $\alpha \in \Sio$ such that the section $B_\alpha$ 
is uncountable.

\hs Recall that by a Theorem of Lusin and Novikov, see \cite[page 123]{Kechris94}, if for
 all  $\alpha \in \Sio$, the section $B_\alpha$ of the Borel set $B$ was   countable,  
 then  $\mathrm{PROJ}_{\Sio}(B)$ would be  a Borel subset of $\Sio$. 

\hs Thus there exists at least  one $\alpha \in \Sio$ such that  $B_\alpha$  is uncountable.
In fact we have not only one  $\alpha$ such that $B_\alpha$  is uncountable.  

\hs For $\alpha \in \Sio$ we have 
$\{\alpha\} \times B_\alpha = B \cap [ \{\alpha\} \times X^\om ]$. 
But $\{\alpha\} \times X^\om$ is a closed 
hence Borel subset of $\Sio \times X^\om$ thus $\{\alpha\} \times B_\alpha$ 
is Borel as intersection of two Borel sets.

\hs If there was only one $\alpha \in \Sio$ such that  $B_\alpha$  is uncountable, then 
$C=\{\alpha\}\times B_\alpha$ would be Borel so $D=B - C$ would be borel 
because the class of Borel sets is closed 
under boolean operations. 
\nl But all sections of $D$ would be countable thus 
 $\mathrm{PROJ}_{\Sio}(D)$ would be Borel by Lusin and Novikov's Theorem. 
Then $\mathrm{PROJ}_{\Sio}(B)= \{\alpha\}\cup \mathrm{PROJ}_{\Sio}(D)$ 
would be also Borel as union of two Borel sets, and this would lead to a
contradiction. 

\hs In a similar manner we can prove that  the set $U=\{\alpha \in \Sio \mid  B_\alpha  
\mbox{ is uncountable } \}$ is uncountable, otherwise 
$U=\{\alpha_0, \alpha_1, \ldots \alpha_n, \ldots \}$ would be Borel as the countable union 
of the closed sets $\{\alpha_i\}$, $i\geq 0$. 
\nl For each $n\geq 0$ the set $\{\alpha_n\}\times B_{\alpha_n}$ would be Borel,  
and 
$C=\cup_{n \in \om}\{\alpha_n\}\times B_{\alpha_n}$ would be Borel as a 
countable union of Borel sets. 
So $D=B - C$ would be borel too. 
\nl But all sections of $D$ would be countable thus 
 $\mathrm{PROJ}_{\Sio}(D)$ would be Borel by Lusin and Novikov's Theorem. 
 Then $\mathrm{PROJ}_{\Sio}(B)= U \cup \mathrm{PROJ}_{\Sio}(D)$ would be 
also Borel as union of two Borel sets, and this would lead to a
contradiction. 

\hs So we have proved that the set 
$\{ \alpha \in \Sio \mid B_\alpha \mbox{ is uncountable } \}$
 is uncountable. 

\hs On the other hand we know from 
another Theorem of Descriptive Set Theory that the set 
$\{ \alpha \in \Sio \mid B_\alpha \mbox{ is countable } \}$ is a ${\bf \Pi^1_1}$-subset of 
$\Sio$, see \cite[page 123]{Kechris94}.  
 Thus its complement $\{ \alpha \in \Sio \mid B_\alpha \mbox{ is uncountable } \}$
is analytic. 
But by Suslin's Theorem an analytic subset of   $\Sio$ is either countable 
or has cardinality $2^{\aleph_0}$, \cite[p. 88]{Kechris94}. Therefore the set 
$\{ \alpha \in \Sio \mid B_\alpha \mbox{ is uncountable } \}$
has cardinality $2^{\aleph_0}$.

\hs Recall now  that we have already seen that, for each $\alpha \in \Sio$, the set 
 $\{\alpha\} \times B_\alpha$ is Borel. 
Thus  $B_\alpha$ 
itself is Borel and  by Suslin's Theorem $B_\alpha$ is either countable or has cardinality $2^{\aleph_0}$.
From this we deduce that 
$\{ \alpha \in \Sio \mid B_\alpha \mbox{ is uncountable } \} = 
\{ \alpha \in \Sio \mid B_\alpha \mbox{ has cardinality } 2^{\aleph_0} \}$ 
has cardinality $2^{\aleph_0}$. \ep

\hs  This Lemma was  used in   \cite{Fink-Sim}  to
 prove  that analytic but non Borel 
context-free $\om$-languages have a maximum degree of ambiguity. 

\begin{The}[\cite{Fink-Sim}]\label{the-amb}
Let $L(\mathcal{A})$ be a context-free $\om$-language  accepted by a B\"uchi pushdown automaton  $\mathcal{A}$ such that $L(\mathcal{A})$ 
is an analytic but non Borel set. Then the  set of $\om$-words, 
which have $2^{\aleph_0}$ accepting runs by $\mathcal{A}$, has cardinality $2^{\aleph_0}$. 
\end{The}

\noi Reasoning in a very similar way as in the proof of Theorem \ref{the-amb} in \cite{Fink-Sim},  we can now state that 
 analytic but non Borel  tree languages accepted by  B\"uchi tree automata have a maximum degree of ambiguity. 

\hs If  $\Si$ is an alphabet having at least two letters, the topological space  $T_\Si^\om$ is homeomorphic to the topological space $\Si^\om$, so we can first 
state Lemma \ref{lemma} in the following equivalent form. 

\begin{Lem}\label{lemma2}
Let $\Si$ and $K$ be two finite alphabets having at least two letters and 
 $B$ be a Borel subset of 
$T_\Si^\om \times T_K^\om$ such that $\mathrm{PROJ}_{T_\Si^\om }(B)$ is not a Borel subset of $T_\Si^\om$.
Then there are $2^{\aleph_0}$ infinite trees   $t \in T_\Si^\om $ such that the section $B_t$ 
has cardinality $2^{\aleph_0}$.
\end{Lem}

\noi We can now state the following result. 

\begin{The}\label{max-amb}
Let $L(\mathcal{A})\subseteq T_\Si^\om$ be a regular tree language  accepted by a B\"uchi tree automaton  $\mathcal{A}$ such that $L(\mathcal{A})$ 
is an analytic but non Borel set. Then the  set of trees $t\in T_\Si^\om$ 
which have $2^{\aleph_0}$ accepting runs by $\mathcal{A}$, has cardinality $2^{\aleph_0}$. 
\end{The}

\proo Let $\mathcal{A}=(K,\Si,\Delta, q_0, F)$ be a B\"uchi tree automaton accepting a non Borel tree language 
$L(\mathcal{A})\subseteq T_\Si^\om$, and let 
$R \subseteq  T_\Si^\om \times T_K^\om$ be defined by : 
$$R= \{ (t, \rho) \mid t\in T_\Si^\om \mbox{ and }   \rho \in T_K^\om  \mbox{ is an accepting run of } \mathcal{A} \mbox{ on  the tree } t \}.$$
\noi The set $R$ can be seen as a tree language over the product alphabet $\Si \times K$. Then it is easy to see that $R$ is accepted by a 
{\it deterministic} B\"uchi tree automaton. But every  tree language  which is accepted by a {\it deterministic} B\"uchi tree automaton is a 
${\bf \Pi}_2^0$-set, see \cite{Murlak05}. 
Thus the tree language $R$ is a ${\bf \Pi}_2^0$-subset of the space  $T_{(\Si\times K)^\om}$ which is 
identified to  the topological  space  $T_\Si^\om \times T_K^\om$. In particular,  $R$ is a Borel subset of $T_\Si^\om \times T_K^\om$. 
But by definition of $R$ it turns out that $\mathrm{PROJ}_{T_\Si^\om}(R)=L(\mathcal{A})$. Thus  
$\mathrm{PROJ}_{T_\Si^\om}(R)$ is not Borel and  Lemma \ref{lemma2} implies that  
there are $2^{\aleph_0}$ trees  $t \in T_\Si^\om$ such that $R_t$ has cardinality 
$2^{\aleph_0}$. This means that these trees have $2^{\aleph_0}$ 
accepting runs by the B\"uchi  tree automaton $\mathcal{A}$. 
\ep 

\begin{Rem}{\rm The above proof  is no longer valid if we replace ``B\"uchi tree automaton" by ``Muller tree automaton". Indeed if 
$L(\mathcal{A})\subseteq T_\Si^\om$ is a regular tree language  accepted by a Muller tree automaton $\mathcal{A}=(K,\Si,\Delta, q_0, \mathcal{F})$, 
 then the set $R \subseteq  T_\Si^\om \times T_K^\om$  defined by : 
$$R= \{ (t, \rho) \mid t\in T_\Si^\om \mbox{ and }   \rho \in T_K^\om  \mbox{ is an accepting run of } \mathcal{A} \mbox{ on  the tree } t \}.$$
\noi is now accepted by a 
deterministic {\bf Muller}  tree automaton. Thus we can now only say that $R$ is a ${\bf \Pi}_1^1$-set,  and we cannot use the fact that $R$ is Borel, which 
was crucial in the proof of Theorem \ref{max-amb}.  
}
\end{Rem}

\noi In particular, Theorem \ref{max-amb} implies the following important result. 

\begin{Cor}\label{unamb}
Let $L(\mathcal{A})\subseteq T_\Si^\om$ be a regular tree language  accepted by an unambiguous  B\"uchi tree automaton. Then the tree language 
$L(\mathcal{A})$ is a Borel subset of $T_\Si^\om$. 
\end{Cor}

\begin{Rem}{\rm The result given by Corollary \ref{unamb} is weaker than the result given by Theorem \ref{max-amb}. This weaker result  can be proved 
by a simpler argument. We give now this proof which is also interesting. 
}
\end{Rem}

\proo  Let $L(\mathcal{A})\subseteq T_\Si^\om$ be a regular tree language  accepted by an unambiguous  
B\"uchi tree automaton $\mathcal{A}=(K,\Si,\Delta, q_0, F)$. 
Let $R$ be defined as in the proof of Theorem \ref{max-amb} by: 
$$R= \{ (t, \rho) \mid t\in T_\Si^\om \mbox{ and }   \rho \in T_K^\om  \mbox{ is an accepting run of } \mathcal{A} \mbox{ on  the tree } t \}.$$
\noi The set $R$ is accepted by a 
{\it deterministic} B\"uchi tree automaton so it is a  ${\bf \Pi}_2^0$-subset of the space  $T_{(\Si\times K)^\om}$. 
\nl Consider now the projection $\mathrm{PROJ}_{T_\Si^\om} : ~ T_\Si^\om \times T_K^\om \ra T_\Si^\om$ defined by 
$\mathrm{PROJ}_{T_\Si^\om}(t, \rho) = t$ for all $(t, \rho) \in  T_\Si^\om \times T_K^\om$. This projection is a continuous function 
and it is {\it injective} on the Borel set $R$ because the automaton $\mathcal{A}$ is unambiguous. 
By a Theorem of Lusin and Souslin, see \cite[Theorem 15.1 page 89]{Kechris94}, the injective 
image of $R$ by the continuous function $\mathrm{PROJ}_{T_\Si^\om}$ is then Borel. 
Thus the tree language $L(\mathcal{A})=\mathrm{PROJ}_{T_\Si^\om}(R)$ is a Borel subset of $T_\Si^\om$. 
\ep 

\begin{Rem}{\rm The above result given by Corollary \ref{unamb}
is of course false in the case of Muller automata because we already know an example of non Borel regular tree language 
accepted by a {\it deterministic hence unambiguous} Muller tree automaton. By Lemma \ref{det}, the tree language 
$\mathcal{L}^-=\fa \mathrm{Path}( \{0, 1\}^\om - (0^\star.1)^\om )$ is a ${\bf \Pi}^1_1$-complete set 
accepted by a deterministic Muller tree automaton. 
}
\end{Rem}

\section{Game tree languages}

\noi Game tree languages are particular recognizable tree languages which are defined by the use of parity games. So we now recall the definition 
of these games, as introduced in  \cite{ArnoldNiwinski08, ADMN}. 

\hs A parity game is a game with perfect information between two players named Eve and Adam, as in \cite{ArnoldNiwinski08, ADMN}.  
\nl  The game is defined by 
a tuple $G=(V_\exists, V_\fa, \mathrm{Move}, p_0,  \mathrm{rank})$. 
The sets  $V_\exists$ and  $V_\fa$ are disjoint sets of positions of Eve and Adam, respectively. We denote 
$V = V_\exists \cup V_\fa$ the set of positions. 
The relation  $\mathrm{Move} \subseteq V \times V$ is the relation of possible moves. The initial position in a play is $p_0 \in V$. The ranking function is 
$ \mathrm{rank} : V \ra \om$ and the number of values taken by this function is finite.  
\nl At the beginning of a play there is a token at the initial position $p_0$ where the play starts. 
The players move the token according to the relation $ \mathrm{Move}$, always to a successor of the current position. 
The move is done by Eve if the current position is an element of $V_\exists$, otherwise Adam moves the token. 
This way  the two players  form a path in the graph $(V,  \mathrm{Move})$. If at some moment a player cannot move then she or he looses. Otherwise the 
two players construct  an infinite path in the graph, $v_0, v_1, v_2, \ldots$ 
In this case Eve wins the play if $\lim \sup _{n\ra \infty}  \mathrm{rank}(v_n)$ is even, otherwise 
Adam wins the play. 
\nl Eve (respectively, Adam) wins the game $G$ if she (respectively, he) has a winning strategy. It is well known that  parity games are determined, 
i. e.,  that one of the players has a winning strategy. Moreover any position is winning for one of the players and she or he has a {\it positional} strategy from 
this position, see \cite{2001automata} for more details. 

\hs We now recall the definition of game languages $W_{(\iota, \kappa)}$. 
\nl  A Mostowski-Rabin index is a pair  $(\iota, \kappa)$, where $\iota \in \{0, 1\}$ and $\iota \leq \kappa < \om$. For such an index, we define the alphabet 
$\Si_{(\iota, \kappa)} = \{ \exists, \fa \} \times \{ \iota, \ldots, \kappa\}$. 
\nl For a letter $a\in \Si_{(\iota, \kappa)}$ we denote $a=(a_1, a_2)$, where $a_1 \in \{ \exists, \fa \}$ and $a_2 \in \{ \iota, \ldots, \kappa\}$. 
\nl  For each tree $t \in T_{\Si_{(\iota, \kappa)}}^\om$ 
we associate a parity game $G(t)=(V_\exists, V_\fa,  \mathrm{Move}, p_0,  \mathrm{rank})$, where 

\begin{itemize}
\ite  $V_\exists=\{ v \in \{l, r\}^\star \mid t(v)_1=\exists \}$,
\ite $V_\fa=\{ v \in \{l, r\}^\star \mid t(v)_1=\fa \}$,
\ite  $ \mathrm{Move} = \{ (w, wi) \mid w \in \{l, r\}^\star \mbox{ and } i \in \{l, r\}\}$,
\ite  $p_0 = \lambda$ is the root of the tree, 
\ite  $ \mathrm{rank}(v)=t(v)_2$, for each $v\in \{l, r\}^\star $.
\end{itemize}

\noi The set      $W_{(\iota, \kappa)} \subseteq T_{\Si_{(\iota, \kappa)}}^\om$   is the set of infinite binary trees $t$ labelled in the alphabet 
$\Si_{(\iota, \kappa)}$ such that Eve wins the associated game $G(t)$. 

\hs The recognizable  tree language $W_{(\iota, \kappa)}$ is accepted by an alternating parity tree automaton of index $(\iota, \kappa)$. This notion 
will be useful in the sequel so we recall it now, as presented in  \cite{ADMN}. 

\begin{Deff} 
An alternating parity tree automaton is a tuple $\mathcal{A}=(\Si, Q_\exists,$ $Q_\fa, q_0, \delta,  \mathrm{rank})$, 
where the set of states $Q$ is partitioned in $Q_\exists$ 
and  $Q_\fa$. The set $Q_\exists$ is the set of existential states and  the set $Q_\fa$ is the set of universal states. The transition relation is 
$\delta \subseteq Q  \times \Si \times \{l, r, \lambda\}  \times Q$ and $ \mathrm{rank}:Q \ra \om$ is the rank function. 
A tree $t \in T_{\Si}^\om$ is accepted by 
the automaton $\mathcal{A}$ iff Eve has a winning strategy in the parity game 
$(Q_\exists \times \{l, r\}^\star, Q_\fa  \times \{l, r\}^\star, (q_0, \lambda),  \mathrm{Move}, 
\Omega)$, where $ \mathrm{Move}=\{ ((p,v), (q,vd)) \mid v \in  \mathrm{dom}(t), ~~ (p, t(v),$ $d, q) \in \delta \}$ and $\Omega(q, v) =  \mathrm{rank}(q)$. 
\end{Deff} 

\noi Notice that it can be assumed without lost of generality
that $ \mathrm{min} ~ \mathrm{rank}(Q)$ is equal to $0$ or $1$. The pair $( \mathrm{min}~ \mathrm{rank}(Q),  \mathrm{max} ~ \mathrm{rank}(Q))$ 
is called the Mostowski-Rabin  index of the automaton. 
\nl It follows from \cite{Rabin69} that any alternating parity tree automaton can be simulated by a non deterministic Muller automaton, 
see also \cite{2001automata}.

\hs There is a usual partial order on Mostowski-Rabin  indices: $(\iota, \kappa) \sqsubseteq (\iota', \kappa') $ 
if either $\iota' \leq \iota$ and $\kappa \leq \kappa'$ (i.e. $\{ \iota, \ldots, \kappa\} \subseteq \{ \iota', \ldots, \kappa'\}$),  
or $ \iota=0,  \iota'=1$, and $\kappa+2 \leq \kappa'$ (i.e. $\{ \iota+2, \ldots, \kappa+2\} \subseteq \{ \iota', \ldots, \kappa'\}$). 
\nl The indices $(1, n)$ and $(0, n-1)$ are called dual and $\overline{(\iota, \kappa)}$ denotes the index dual to $(\iota, \kappa)$.

\hs It is easy to see 
that each tree language $W_{(\iota, \kappa)}$ is accepted by an alternating parity tree automaton of index $(\iota, \kappa)$.

\hs Moreover the set $W_{(\iota, \kappa)}$ is in some sense of the greatest possible topological complexity among tree languages accepted by 
alternating parity tree automata of index $(\iota, \kappa)$. This is expressed by the following lemma. 

\begin{Lem}[ see \cite{ADMN} ]\label{lemma-index}
If a set of trees $T$ is recognized by an  alternating parity tree automaton of index $(\iota, \kappa)$, then $T \leq_W W_{(\iota, \kappa)}$. 
\end{Lem}

\noi In order to use this result to get a lower bound on the  topological complexity of the game tree languages $W_{(\iota, \kappa)}$, we first 
construct some alternating parity tree automata accepting the tree languages  $\mathcal{L}$ and  $\mathcal{L}^-$  defined in the preceding section. 

\begin{Lem}
The tree language $\mathcal{L}$ is accepted by an alternating parity tree automaton of index $(1, 2)$. 
\end{Lem}

\proo Recall that $\mathcal{L}= \exists \mathrm{Path}( \mathcal{R} )$, where $\mathcal{R}=(0^\star.1)^\om$. 
\nl The tree language $\mathcal{L}$ is then accepted by the 
alternating parity tree automaton $\mathcal{A}=(\Si, Q_\exists,$ $Q_\fa, q_0, \delta,  \mathrm{rank})$, where 
\nl $\Si=\{0,1\}$, 
\nl $Q_\exists=Q=\{q_0, q_1\}$, 
\nl $Q_\fa=\emptyset$, 
\nl $\delta=\{(q, 1, d, q_1), (q, 0, d, q_0) \mid q\in Q \mbox{ and } d\in\{l, r\} \}$,
\nl $ \mathrm{rank}(q_0)=1$ and $ \mathrm{rank}(q_1)=2$. 
\ep 

\hs Notice that  in the above automaton $\mathcal{A}$ all states are existential. 

\begin{Lem}
The tree language $\mathcal{L}^-$ is accepted by an alternating parity tree automaton of index $(0, 1)$. 
\end{Lem}

\proo Recall that $\mathcal{L}^-=T_\Si^\om - \mathcal{L}=\fa  \mathrm{Path}( \{0, 1\}^\om - (0^\star.1)^\om )$. 
\nl The tree language $\mathcal{L}^-$ is then accepted by the 
alternating parity tree automaton $\mathcal{A}'=(\Si, Q'_\exists,$ $Q'_\fa, q'_0, \delta',  \mathrm{rank}')$, where 
\nl $\Si=\{0,1\}$, 
\nl $Q'_\exists=\emptyset$, 
\nl $Q'_\fa=Q'=\{q'_0, q'_1\}$, 
\nl $\delta'=\{(q', 1, d, q'_1), (q', 0, d, q'_0) \mid q'\in Q' \mbox{ and } d\in\{l, r\} \}$,
\nl $\mathrm{rank}'(q'_0)=0$ and $ \mathrm{rank}'(q'_1)=1$. 

\ep 

\hs Notice that  in the above automaton $\mathcal{A}'$ all states are universal. 

\begin{Rem}
The ${\bf \Si}^1_1$-complete tree language $\mathcal{L}$ is accepted by an alternating parity tree automaton of index $(1, 2)$ and the 
${\bf \Pi}^1_1$-complete tree language $\mathcal{L}^-$ is accepted by an alternating parity tree automaton of index $(0, 1)$. In fact for 
every tree language $T$ accepted by an  alternating parity tree automaton of index $(1, 2)$ (respectively, $(0, 1)$) it holds that 
$T$ is in the class ${\bf \Si}^1_1$ (respectively,  ${\bf \Pi}^1_1$), see \cite[Theorem 3.6]{ADMN}. 
\end{Rem}

\noi Recall now the definition of the $D_{\om}({\bf \Si}^1_1)$-complete    tree language  $\mathcal{L}_1$. 

\hs  $\mathcal{L}_1=\{ t \in T_{\{0, 1\}}^\om  \mid  \exists n\geq 0 ~~t_{l^n\cdot r} \in \mathcal{L} \mbox{ and  min} 
\{ n\geq 0 \mid t_{l^n\cdot r} \in \mathcal{L} \} 
\mbox{ is odd }\}$. 

\hs We can now state the following result. 

\begin{Lem}
The tree language $\mathcal{L}_1$ is accepted by an alternating parity tree automaton of index $(0, 2)$. 
\end{Lem}

\proo Let, as in the proofs of the two previous lemmas, $\mathcal{A}=(\Si, Q_\exists,$ $Q_\fa, q_0, \delta,  \mathrm{rank})$ be an 
alternating parity tree automaton of index $(1, 2)$ accepting the tree language $\mathcal{L}=\exists  \mathrm{Path}( \mathcal{R} )$, and 
$\mathcal{A}'=(\Si, Q'_\exists,$ $Q'_\fa, q'_0, \delta',  \mathrm{rank}')$  be an 
alternating parity tree automaton of index $(0, 1)$ accepting the tree language $\mathcal{L}^-$. We assume that 
$Q \cap Q'=\emptyset$, where $Q=Q_\exists \cup Q_\fa=Q_\exists$ and $Q'=Q'_\exists \cup Q'_\fa=Q'_\fa$. 
\nl It is then easy to see  that the tree language  $\mathcal{L}_1$
is accepted by the 
alternating parity tree automaton $\mathcal{A}^1=(\Si, Q^1_\exists,$ $Q^1_\fa, q^1_0, \delta^1,  \mathrm{rank}^1)$, where 
\nl $\Si=\{0,1\}$, 
\nl $Q^1_\exists = Q_\exists \cup  Q'_\exists \cup \{q_\exists\}= Q_\exists \cup  \{q_\exists\}$, 
\nl $Q^1_\fa = Q_\fa   \cup Q'_\fa  \cup \{q^1_0, q^1_1\}= Q'_\fa  \cup \{q^1_0, q^1_1\}$, 
\nl $\delta^1= \delta \cup \delta' \cup            \{(q^1_0, a, l, q_\exists),  (q^1_0, a ,r ,q'_0), (q_\exists, a, r,  q_0),  
(q_\exists, a, \lambda,  q^1_1) ,   (q^1_1, a, r, q'_0),  (q^1_1, a, l, q^1_0) 
  \mid a\in \{0, 1\} \}$,
\nl  $ \mathrm{rank}^1(q)= \mathrm{rank}(q)$ for $q\in Q$, 
\nl  $ \mathrm{rank}^1(q')= \mathrm{rank}'(q')$ for $q'\in Q'$, 
\nl  $ \mathrm{rank}^1(q^1_0)=0$, $ \mathrm{rank}^1(q^1_1)=1$. 

\ep 

\hs Notice that in the above construction of the alternating automaton $\mathcal{A}^1$ the universal states $q^1_0, q^1_1$ and the existential state 
$q_\exists$ are used to choose, when reading a tree $t\in \mathcal{L}_1$,  the least integer $n$ such that $t_{l^n\cdot r} \in \mathcal{L}$ and to check 
that this integer is really the least (and odd) one with this property. 
\nl In a very similar manner, for each integer $n\geq 1$, 
 we can define an alternating parity tree automaton $\mathcal{A}^n$ of index $(0, 2)$ accepting the language $\mathcal{L}_n$. 
The complete description would be tedious  but the idea is that now the additional universal or existential states not in $Q \cup Q'$ 
are used to choose, for a given tree    $t\in \mathcal{L}_n$,   
  the least ordinal $\alpha= \om^{n-1}\cdot a_{n-1} + \om^{n-2}\cdot a_{n-2} + \ldots + \om\cdot a_1 + a_0$ 
such that $t_{l^{a_{n-1}}\cdot r\cdot l^{a_{n-2}}\cdot r \cdots l^{a_0}\cdot r}$ is in  $\mathcal{L}$ and to check that $\alpha$ is odd and that for any 
smaller ordinal $\beta= \mathrm{Ord}(b_{n-1}, b_{n-2}, \ldots , b_0) < \alpha$,  the tree  
 $t_{l^{b_{n-1}}\cdot r \cdot l^{b_{n-2}}\cdot r \cdots l^{b_0}\cdot r}$ is not in  $\mathcal{L}$. 

\hs We can then state the following result. 

\begin{Pro}\label{Ln-alt}
For each integer $n\geq 1$, the tree language $\mathcal{L}_n$ is accepted by an alternating parity tree automaton of index $(0, 2)$. 
\end{Pro}

\noi We can now infer from Theorem \ref{Ln}, Proposition \ref{Ln-alt},   and Lemma \ref{lemma-index}, the following result. 

\begin{The}\label{main-theor}
For each integer $n \geq 1$, the $D_{\om^n}({\bf \Si}^1_1)$-complete tree language $\mathcal{L}_n$ is Wadge reducible to the game tree language
$W_{(0,2)}$, i.e. ~~$\mathcal{L}_n \leq_W W_{(0,2)}$. In particular the language  $W_{(0,2)}$ is not in any class $D_{\alpha}({\bf \Si}^1_1)$ 
for $\alpha < \om^\om$. 
\end{The}

\noi On the other hand, Arnold and Niwinski proved in \cite{ArnoldNiwinski08} that the game tree languages form a hierarchy with regard to the Wadge 
reducibility. 

\begin{The}[\cite{ArnoldNiwinski08}]\label{strict-hierarchy}
For all Mostowski-Rabin indices  $(\iota, \kappa)$ and $(\iota', \kappa')$, it holds that : 
$$(\iota, \kappa)    \sqsubseteq   (\iota', \kappa')  ~~  \mbox{ if and only if } ~~W_{(\iota, \kappa)} \leq_W  W_{(\iota', \kappa')} $$
\end{The}

\noi Then we can state the following result. 

\begin{The}
For each integer $n \geq 1$ and each Mostowski-Rabin index  $(\iota, \kappa)$ 
such that $(0,2) \sqsubseteq(\iota, \kappa)$ or $(\iota, \kappa)=(1,3)=\overline{(0,2) }$, 
 the $D_{\om^n}({\bf \Si}^1_1)$-complete tree language $\mathcal{L}_n$ is Wadge reducible to the game tree language
$W_{(\iota, \kappa)}$, i.e. ~~$\mathcal{L}_n \leq_W W_{(\iota, \kappa)}$. 
In particular the language  $W_{(\iota, \kappa)}$ is not in any class $D_{\alpha}({\bf \Si}^1_1)$ 
for $\alpha < \om^\om$. 
\end{The}

\proo The result follows directly from Theorems \ref{main-theor} and \ref{strict-hierarchy} in the case $(0,2) \sqsubseteq(\iota, \kappa)$. 
What remains is the case of the index $(1,3)$ which is the dual of the index $(0,2)$. But it is proved in 
\cite[Lemma 1]{ArnoldNiwinski08} that $W_{\overline{(\iota, \kappa)}}$ 
coincide with $\overline{W_{(\iota, \kappa)}}= T_{\Si_{(\iota, \kappa)}}^\om - W_{(\iota, \kappa)}$ up to renaming of symbols. On the other hand, 
we know from 
Theorem \ref{main-theor} that for each integer $n \geq 1$, 
the $D_{\om^{n+1}}({\bf \Si}^1_1)$-complete tree language $\mathcal{L}_{n+1}$ is Wadge reducible to the game tree language
$W_{(0,2)}$, i.e. ~~$\mathcal{L}_{n+1} \leq_W W_{(0,2)}$. This is easily seen to be equivalent to 
$T_{\{0, 1\}}^\om - \mathcal{L}_{n +1} \leq_W \overline{W_{(0,2)}}$, i.e. $T_{\{0, 1\}}^\om - \mathcal{L}_{n+1}  \leq_W W_{(1,3)}$. 
But $\mathcal{L}_{n}$ is $D_{\om^{n}}({\bf \Si}^1_1)$-complete and $\mathcal{L}_{n+1}$ is $D_{\om^{n+1}}({\bf \Si}^1_1)$-complete 
so it follows from the properties of the  difference hierarchy of analytic sets that $\mathcal{L}_{n}  \leq_W T_{\{0, 1\}}^\om - \mathcal{L}_{n+1}$ and so 
$\mathcal{L}_{n}  \leq_W W_{(1,3)}$ by transitivity of the relation $\leq_W$.
\ep 

\section{Concluding remarks}

\noi We have got some new results on 
 the topological complexity of  non Borel recognizable tree languages with regard to the difference hierarchy of analytic sets. In 
particular, we have showed that the game tree language  $W_{(0,2)}$ is not in any class $D_{\alpha}({\bf \Si}^1_1)$ 
for $\alpha < \om^\om$. 
The great challenge in the study of the topological complexity of  recognizable tree languages is to determine the Wadge hierarchy of 
tree languages accepted by {\it non deterministic} Muller or Rabin tree automata. Notice that the case of {\it  deterministic} Muller or Rabin tree automata
have been solved recently by Murlak, \cite{Murlak-LMCS}.  

\hs  It would be interesting to locate in a more precise way the game tree languages with regard to the difference hierarchy of analytic sets. 
We already know that  $W_{(0,2)}$ is not in any class $D_{\alpha}({\bf \Si}^1_1)$ for $\alpha < \om^\om$.  Is there an ordinal $\alpha$ such that 
 $W_{(0,2)}$ is in $D_{\alpha}({\bf \Si}^1_1)$ and then what is the smallest such ordinal $\alpha$?  The same question may be asked for the other 
game tree 
languages $W_{(\iota, \kappa)}$. 
 On the other hand, there are some sets in the class ${\bf \Delta}^1_2$ which does not belong to the 
$\sigma$-algebra generated by the analytic sets, see \cite[Exercise 37.8]{Kechris94}. Could we expect that $W_{(0,2)}$ or another game tree language 
 $W_{(\iota, \kappa)}$  is such an example?

\hs  {\bf  Acknowledgements.}  We thank the anonymous referees for their  very helpful  comments which have led to a great improvement  of our paper.

\end{document}